\documentclass{amsart}
\usepackage{amsaddr}
\usepackage{hyperref}
\usepackage{etoolbox}
\usepackage{lipsum}

\usepackage{graphicx}
\usepackage{tikz}
\usepackage{tikz-cd}
\usepackage{verbatim}
\usepackage{setspace}
\usepackage{mathrsfs}
\usepackage{hyperref}
\usepackage{academicons}
\usepackage{xcolor}

\newtheorem{theorem}{Theorem}[section]
\newtheorem{lemma}[theorem]{Lemma}

\theoremstyle{definition}
\newtheorem{definition}[theorem]{Definition}

\newtheorem{examples}[theorem]{Examples}

\theoremstyle{remark}

\theoremstyle{proposition}
\newtheorem{proposition}[theorem]{Proposition}

\theoremstyle{corollary}
\newtheorem{corollary}[theorem]{Corollary}

\numberwithin{equation}{section}

\DeclareMathOperator{\Hom}{Hom} 
 \DeclareMathOperator{\ann}{ann}

\DeclareMathOperator{\E}{E}

\newcommand{\abs}[1]{\lvert#1\rvert}

\newcommand*\xbar[1]{%
\hbox{%
\vbox{%
\hrule height 0.5pt 
\kern0.5ex
\hbox{%
\kern-0.1em
\ensuremath{#1}%
\kern-0.1em
}%
}%
}%
}
\newcommand\restr[2]{{
\left.\kern-\nulldelimiterspace 
#1 
\right|_{#2} 
}}

\newcommand{\orcid}[1]{\href{https://orcid.org/#1}{\textcolor[HTML]{A6CE39}{\aiOrcid}}}

\begin{document}
\title{On $R$-Coneat Injective Modules and Generalizations}

\author{Mohanad Farhan Hamid} 
\address{Department of Applied Siences,
\\University of Technology-Iraq, Baghdad, Iraq
\\E-mails: {mohanad.f.hamid@uotechnology.edu.iq
\newline mohanadfhamid@yahoo.com}}

 \curraddr{}

\thanks{}





\maketitle

\section*{\bf{Abstract}}
\noindent Both the classes of $R$-coneat injective modules and its superclass, pure Baer injective modules, are shown to be preenveloping. The former class is contained in another one, namely, self coneat injectives, i.e. modules $M$ for which every map $f$ from a coneat left ideal of $R$ into $M$, whose kernel contains the annihilator of some element in $M$, is induced by a homomorphism $R \rightarrow M$. Certain types of rings are characterized by properties of the above modules. For instance, a commutative ring $R$ is von Neuman regular if and only if all self coneat injective $R$-modules are quasi injective.

\bigskip
\noindent \textbf{\keywordsname.} coneat injective module, copure Baer injective module, right SF-ring.

\section{\bf{Introduction}} The ring $R$ is always associative with identity. Unless otherwise stated, modules are left unital $R$-modules. Homomorphisms and maps are all $R$-homomorphisms. 
A sequence of left $R$-modules: 
\begin{equation}\label{seq}
0 \rightarrow A \hookrightarrow B \rightarrow B/A \rightarrow 0
\end{equation}
is \emph{pure} (respectively, \emph{coneat}) exact if the sequence $0 \rightarrow U\otimes A \rightarrow U\otimes B \rightarrow U\otimes B/A \rightarrow 0$ remains exact for all (simple) right $R$-modules $U$. Therefore, pure exact sequences are coneat exact but not convesely. For a nonempty class of modules $\mathcal Q$, the sequence \ref{seq} is called $\mathcal Q$-\emph{copure} exact \cite[section 38]{wis} if every module $Q \in \mathcal Q$ is injective with respect to \ref{seq}, i.e. every map $A \rightarrow Q$ has an extension to a map $B \rightarrow Q$.
A module is \emph{injective} (respectively, \emph{pure} injective, \emph{coneat} injective) if it is injective with respect to all exact (respectively, pure exact, coneat exact) sequences  \ref{seq}.  Baer condition states that the module remains injective if we restrict to the case when $B$ in the above sequence is the ring $R$. No `Baer condition' is known for pure injectivity or coneat injectivity. Therefore, we have \emph{pure Baer injective} and $R$-\emph{coneat injective} modules, i.e. modules injective with respect to all pure and coneat exact sequences \ref{seq} with $B=R$, respectively \cite{nada, cinj}. Also, $\mathcal Q$-\emph{copure Baer injective} modules are those injective with respect to all $\mathcal Q$-copure exact sequences \ref{seq} with $B=R$ \cite{cpbi}. By \cite[p. 290]{wis} and \cite[Theorem 2.3]{cinj}, pure exact and coneat exact sequences are $\mathcal Q$-copure exact by taking $\mathcal Q=$ the class of pure injective or coneat injective modules, accordingly. Hence, pure Baer and $R$-coneat injectivity are examples of $\mathcal Q$-copure Baer injectivity.

The module $M$ is \emph{quasi injective} \cite{jwong} if it is injecitve with respect to every sequence \ref{seq} with $B =M$. L. Fuchs \cite{fuchs} gave a kind of Baer condition for quasi injectivity. He proved that for a module $M$  to be quasi injecitve, it is necessary and sufficient that every map $L \rightarrow M$ from a left ideal $L$ of $R$, the kernel of which is in $\Omega (M)$ (equivalentley, $\xbar{\Omega} (M)$), can be extended to a map $R \rightarrow M$. Here, the set $\Omega (M)$ consists of all left ideals $L$ of $R$ such that $L \supseteq \ann (m)$ for some $m \in M$, and $\xbar{\Omega} (M)$ is the filter generated by $\Omega (M)$. This motivated the study of absolutely self pure, absolutely self neat and principally self injective modules \cite{absp, absn, psi}. Analogously, we study here \emph{self coneat injective} modules by taking maps $L \rightarrow M$, from coneat left ideal $L$ of $R$, whose kernels are in $\Omega (M)$, extendable to $R \rightarrow M$.

A class of modules $\mathcal{P}$ is \emph {preenveloping} if every module $M$ has a $\mathcal{P}$-\emph {preenvelope}, i.e. a map $\varphi : M \rightarrow P$ where $P \in \mathcal{P}$ such that for every map $\psi:M \rightarrow P'$, where $P' \in \mathcal{P}$, there is a map $f: P \rightarrow P'$ with $f \circ \varphi =\psi$. A $\mathcal{P}$-preenvelope $\varphi : M \rightarrow P$ is a $\mathcal{P}$-\emph{envelope} if for $P' = P$ and $\psi = \varphi$ the only such maps $f$ are automorphisms of $P$. The class $\mathcal P$ is called \emph{enveloping} if every module has a $\mathcal P$-envelope. If $\mathcal P$ contains the injective modules then any $\mathcal P$-preenvelope is a monomorphism. In this case, we say that the module $P \in \mathcal P$, rather than the map $\varphi : M \rightarrow P$, is a $\mathcal P$-(pre)envelope of $M$. For details about (pre)envelopes in general the reader is referred to \cite{eandj}. The classes of injective, pure injective, and coneat injective modules are all enveloping. In \cite{cpbi} it is proved that every module is embeddable as a $\mathcal Q$-copure submodule of a $\mathcal Q$-copure Baer injective module. Here we will show that such embedding is actually a $\mathcal Q$-copure Baer injective preenvelope. Thus, both the classes of $R$-coneat injective modules and pure Baer injective modules are preenveloping.
\noindent 
\section{\bf $R$-Coneat Injective Preenvelopes}
\noindent 
In \cite[Theorem 2, Lemma 2]{cpbi}, it is proved that, for any nonempty class of modules $\mathcal Q$, every module can be embedded as a $\mathcal Q$-copure submodule of a $\mathcal Q$-copure Baer injective module. Hence, every module can be embedded as a coneat submodule of  an $R$-coneat injective module. In the following lemma, we modify the proof of \cite[Lemma 2]{cpbi} to show that the above embedding is actually a preenvelope.
\begin{lemma} \label{embeddinglemma} Every module has an $R$-coneat injective preenvelope.
\begin{proof} Given a module $M$, we will embed $M$ in an $R$-coneat injective preenvelope. Consider the coneat left ideals $C$ of $R$ and the set $\mathscr F$ of all maps $f:C \rightarrow M$. For every $f \in \mathscr F$ there is a pushout $P=(R \oplus M)/K$, where $K=\{(c,-f(c)), c\in C\}$, and a map $g:R \rightarrow P$ with $\restr{g}{C} = f$.

For each map $\alpha : M \rightarrow N$ into an $R$-coneat injective module $N$ take the diagram:
$$\begin{tikzcd}
C \arrow[hook]{r}\arrow{d}{f}
&R \arrow{d}{}
\\
M \arrow{r}{h} \arrow{d}{\alpha}& P \arrow[dashed]{ld}{\gamma}\\
N
\end{tikzcd}$$
Since $N$ is $R$-coneat injective, there is a $\beta :R \rightarrow N$ extending $\alpha f$. Define $\gamma : P \rightarrow N$ as: $\gamma ((r,m) + K) = \beta (r) + \alpha (m)$. It is easy to check that $\gamma$ is a  homomorphism extending $\alpha$.

The module $P$ might not be $R$-coneat injective, so put $M_0 = M$, $M_1 = P$, $f_0 = f$, $h_0 = h$, $\alpha _0 = \alpha$ and $\gamma _0 = \gamma$ and repeat the above process with $M$ replaced by $M_1$ to give $M_2$ and $M_0 \subseteq M_1 \subseteq M_2$ and the commutative diagram:
$$\begin{tikzcd}
 C \arrow[hook]{r}\arrow{d}{f_1 = h_0 f_0}
&R \arrow{d}{}
\\
M_1 \arrow{r}{h_1} \arrow{d}{\gamma_1}& M_2 \arrow[dashed]{ld}{\gamma_2}\\
N
\end{tikzcd}$$
Continue in this manner to get a sequence $M_0 \subseteq \cdots \subseteq M_n \subseteq M_{n+1} \subseteq \cdots$, for all $n \in \mathbb{N}$. Put $M_\omega = \bigcup M_n$. Now, for each nonlimit ordinal redo the above process. If we reach a limit ordinal, say $\lambda$, put $M_\lambda = \bigcup \{M_\rho, \rho \, \textless \, \lambda \}$. Let $\sigma$ be the smallest ordinal with a larger cardinality than that of the ring $R$, i.e. $\abs{\sigma} =\abs{R}^+$ (the successor cardinal of $\abs{R}$). For each $\rho \, \textless \, \sigma$, we have $\abs{\rho} \thinspace \textless \thinspace \abs{\sigma}$, $\sigma$ is a limit ordinal and $M_\sigma = \bigcup \{M_\rho, \rho \thinspace \textless \thinspace \sigma \}$.  Now, $M_\sigma$ is our $R$-coneat injective module. To see this, let $C$ be a coneat left ideal of $R$ and $f: C \rightarrow M_\sigma$ any map.  For each $c \in C$, let $\bar c$ be the smallest ordinal such that $f(c) \in M_{\bar c}$. Then $\bar c \thinspace \textless \thinspace \sigma$ and $\abs{\bar c} \thinspace \textless \thinspace \abs{\sigma} =\abs{R}^+$. Hence $\abs{\bar c} \leq \abs{R}$. Put $p = \sup \{\bar c, c \in R\}$. As each $\abs{\bar c} \leq \abs{R}$, we must have $\abs{p} \leq \abs{R} \, \textless \,\, \abs{\sigma}$. Hence, $p \thinspace \textless \thinspace \sigma$. Since $\sigma$ is a limit ordinal, we have $p+1 \thinspace \textless \thinspace \sigma$. Therefore, for each $c \in R$, $c \in M_{\bar c} \subseteq M_p \subseteq M_{p+1} \subseteq M_\sigma$. So, $f(C) \subseteq M_p$. Moreover, the map $f:C \rightarrow M_p$ can be extended to a map $g:R \rightarrow M_{p+1}$ with $\restr{g}{C} = f$. View $g$ now as a map $R \rightarrow M_\sigma$. The map $\gamma : M_\sigma \rightarrow N$ is defined as follows. For every $x \in M_\sigma$ there is an $M_\rho$ with $x \in M_\rho$. So put $\gamma(x) = \gamma_\rho (x)$.
\end{proof}
\end{lemma}

Replacing the word `coneat' by `$\mathcal Q$-copure' in the proof of above Lemma, we get:
\begin{lemma} \label{embedcopure} Every module has a $\mathcal Q$-copure Baer injective preenvelope.
\end{lemma}

Now this Lemma together with \cite[Lemmas 3 and 4]{cpbi} are all we need to show:
\begin{theorem} For a given nonempty class of modules $\mathcal{Q}$, every module can be embedded as a $\mathcal{Q}$-copure submodule in a $\mathcal{Q}$-copure Baer injective preenvelope.
\end{theorem}
\begin{corollary} 
\begin{enumerate}
       \item Every module can be embedded as a coneat submodule of an     
              $R$-coneat injective preenvelope.
       \item Every module can be embedded as a pure submodule of a     
              pure Baer injective preenvelope.
\end{enumerate}
\end{corollary}
\section{\bf Self Coneat Injectivity}
\begin{definition} An $R$-module $M$ is called \emph {self coneat injective} if any map $C \rightarrow M$ from a coneat left  ideal $C$ of $R$, the kernel of which contains $\ann (m)$ for some $m \in M$, can be extended to a map $R \rightarrow M$.
\end{definition}

The argument that a direct summand of a self coneat injective module is self coneat injective is standard.

As $\Omega (R)$ consists of all left ideals of $R$, we have any module containing the ring $R$ is self coneat injective if and only if it is  $R$-coneat injective. In particular, any preenvelope of $R$  is self coneat injective if and only if it is $R$-coneat injective. 
The following Theorem characterizes self coneat  injective modules.
\begin{theorem} An $R$-module $M$ is self coneat injective if and only if for any coneat left ideal $C$ of $R$ and any $m \in M$, every map $Cm \rightarrow M$ has an extension to $Rm \rightarrow M$.
\begin{proof} ($\Rightarrow$) Let $g \in \Hom (Cm,M)$ and define $f: C \rightarrow M$ by $x \mapsto g(xm)$. 
For any $x \in \ann (m)$ we have $f(x) = g(xm) =0$ and hence $x \in \ker f$. Therefore, by assumption, there is an $\bar f \in \Hom (Rm,M)$ extending $f$. 
 The map $Rm \rightarrow M$ given by $m \mapsto \bar f (1)$ is easily seen to be an extension of $g$. ($\Leftarrow$) Let now $f: C \rightarrow M$ be a map from a coneat left ideal of $R$ with kernel containing $\ann (m)$ for some $m \in M$. If $xm=0$ for some $x \in C$ then $x \in \ann (m) \subseteq \ker f$. Therefore, the map $g:Cm \rightarrow M$ given by $xm \mapsto f(x)$ is well defined. By assumption, $g$ can be extended to a $\bar g : Rm \rightarrow M$. Now the map $R \rightarrow M$, $1 \mapsto \bar g (m)$ extends $f$.
\end{proof}
\end{theorem}
\begin{theorem} \label{hereditary}
The following statements are equivalent for a ring $R$:
\begin{enumerate}
       \item Every coneat left ideal of $R$ is projective.
       \item Homomorphic images of $R$-coneat injective modules are again $R$-coneat injective.
       \item Homomorphic images of $R$-coneat injective modules are self coneat injective.
\end{enumerate}
\begin{proof}
(1) $\Rightarrow$ (2)	Let $M$ be an $R$-coneat injective module and $\pi: M \rightarrow N$ an epimorphism. For each coneat left ideal $C$ of $R$ consider the following diagram:
$$\begin{tikzcd}
R 
&C \arrow[l,hook'] {i}\arrow{d}{f} \\
M\arrow{r}{\pi}&N
\end{tikzcd}$$
Since $C$ is projective, there is a lifting map $\alpha : C \rightarrow M$ for $f$. Since $M$ is $R$-coneat injective, there is an extension $\beta : R \rightarrow M$ of $\alpha$. Now $\pi \alpha$ extends $f$. (2) $\Rightarrow$ (3) is trivial. (3) $\Rightarrow$ (1) Let $C$ be a coneat left ideal and consider a diagram like the above with $M$ being injective and, by assumption, $N$ is self coneat injective. Any map $f$ whose kernel contains the annihilator of some $x \in N$ is induced by a map $\alpha : R \rightarrow N$. By projectivity of $R$ there is a $\beta : R \rightarrow M$ lifting $\alpha$. Now $\pi \beta$ lifts $\alpha$, hence $f$, and $C$ is projective by \cite[Lemma 3.4]{absp}
\end{proof}
\end{theorem}

Recall that  right SF-rings $R$ are characterized by the property that every left ideal of $R$ is coneat \cite[Theorem 3.16]{cinj}. Such rings are characterized using self coneat injectivitey as follows.  
\begin{theorem} \label{Rcinjinj}
The following five statements are equivalent for a ring $R$:
\begin{enumerate}
    \item The ring $R$ is a right SF-ring.
    \item Every left ideal of $R$ is coneat.
    \item Every $R$-coneat injective module is injective.
    \item All $R$-coneat injective modules are quasi injective.
    \item All self coneat injective modules are quasi injective.

If $R$ is commutative then the above are equivalent to:
    \item The ring $R$ is von Neuman regular. 
\end{enumerate}
\begin{proof} (1) $\Leftrightarrow$ (2) and in the commutative case (1) $\Leftrightarrow$ (6) By \cite[Theorem 3.16]{cinj}. (2) $\Leftrightarrow$ (3) By \cite[Theorem 3]{cpbi} since coneat submodules are copure. (2) $\Rightarrow$ (4) and (2) $\Rightarrow$ (5) If all left ideals are coneat then by \cite[Theorem 3.16]{cinj}, (4) and (5) follow. (3) $\Rightarrow$ (2) In particular, all coneat injective modules are injective and by \cite[Theorem 3.16]{cinj} (1) follows. (3) $\Rightarrow$ (4) and (5) $\Rightarrow$ (4) are trivial. (4) $\Rightarrow$ (3) Let $M$ be $R$-coneat injective, then so is $M \oplus P$ where $P$ is a preenvelope of $R$. Now, any map from a coneat left ideal $L$ of $R$ into $M$ has an extension to a map $L \rightarrow M \oplus P$ whose kernel contains the annihilator of $(0,1)$. By assumption $M \oplus P$ is quasi injective.
\end{proof}
\end{theorem}

\begin{examples}
\begin{enumerate}
     \item Any quasi injective module and any ($R$-)coneat injective module is self coneat injective.
     \item Since  $\mathbb Z$ is not an SF-ring there must exist then, by Theorem \ref{Rcinjinj}, a self coneat injective $\mathbb Z$-module that is not quasi injective.
     \item If $R$ is a commutative von Neuman regular but not noetherian ring
then there must exist a quasi injective, hence self coneat injective, $R$-module that is
not injective \cite[Proposition 1]{byrd}, hence not $R$-coneat injective. \item Let $A$ be a self coneat injective $\mathbb Z$-module that is not quasi injective  (example 2 above) and $B$ a quasi injective $R$-module that is not $R$-coneat injective (example 3 above).
The self coneat injective module
$\left( {\begin{array}{c}A  \\
    B \\
  \end{array} } \right)$ over the ring
   $\left( {\begin{array}{cc}
   \mathbb Z & 0 \\
   0 &  R\\
  \end{array} } \right)$ is not quasi injective nor $R$-coneat injective, for otherwise the direct summand $\left( {\begin{array}{c}A  \\
    0 \\
  \end{array} } \right)$ would be quasi injective or the direct summand $\left( {\begin{array}{c}0  \\
    B \\
  \end{array} } \right)$ would be $R$-coneat injective. 
\end{enumerate}
\end{examples}

\begin{theorem} \label{csplit}
The following statements are equivalent for a ring $R$:
\begin{enumerate}
       \item Every coneat exact sequence of left $R$-modules splits.
       \item All $R$-modules are $R$-coneat injective.
       \item All $R$-modules are self coneat injective.
\end{enumerate}
\begin{proof} (1) $\Rightarrow$ (2) $\Rightarrow$ (3) are trivial. (3) $\Rightarrow$ (1) Let $C$ be a coneat left ideal of $R$ and $P$ a preenvelope of $R$. By assumption, $C \oplus P$ is self coneat injective. Therefore, as the kernel of the inclusion $C \rightarrow C \oplus P$ contains the annihilator of $(0,1)$, it must be induced by a map $R \rightarrow C \oplus P$. This map followed by the projection $C \oplus P \rightarrow C$ gives a homomrphism $R \rightarrow C$ extending the identity of $C$.
\end{proof}
\end{theorem}

We know that all pure exact sequences of left $R$-modules are coneat exact. For the converse, we have the following.

\begin{proposition} \label{coneatpure}
The following statements are equivalent for a ring $R$:
\begin{enumerate}
       \item Every coneat exact sequence of left $R$-modules is pure exact.
       \item All pure Baer injective $R$-modules are $R$-coneat injective.
       \item All pure Baer injective $R$-modules are self coneat injective.
       \item All self pure injective $R$-modules (i.e. modules $M$ for which every map $f$ from a pure left ideal of $R$ into $M$ whose kernel contains the annihilator of some element in $M$ is induced by a homomorphism $R \rightarrow M$) are self coneat injective.
\end{enumerate}       
\begin{proof} (1) $\Rightarrow$ (2) $\Rightarrow$ (3) and (1) $\Rightarrow$ (4) $\Rightarrow$ (3) are trivial. (2) $\Rightarrow$ (1) By assumption, any pure Baer injective module $M$ is injective with respect to any coneat exact sequence $0 \rightarrow A \rightarrow B \rightarrow C \rightarrow 0$ of left $R$-modules. Therefore, as pure submodules are copure \cite{cpbi}, the above sequence must be pure exact \cite[Theorem 3]{cpbi}. (3) $\Rightarrow$ (2) Let $M$ be pure Baer injective, therefore, so is $M \oplus P$, where $P$ is a pure Baer injective preenvelope of $R$. By assumption,  $M \oplus P$ is self coneat injective. So for any map $f: C \rightarrow M$, the extension $C \rightarrow M \hookrightarrow M \oplus P$ has kernel containing the annihilator of $(0,1)$ and hence, it is induced by a map $R \rightarrow M \oplus P$. Following this last map by the projection $M \oplus P \rightarrow M$, we get an extension $R \rightarrow M$ of $f$.
\end{proof}
\end{proposition}

By similar arguments to that of \ref{coneatpure}, we get the following three propositions.
For the first proposition, recall that a module $M$ is called \emph{finitely} $R$-injecitve (\emph{absolutely self pure}) \cite{rama} (\cite{absp}) if every map $I \rightarrow M$ from a finitely generated left ideal of $R$ (whose kernel is in $\xbar {\Omega} (M)$) can be extended to $R \rightarrow M$.

\begin{proposition} \label{coneatfg}
The following statements are equivalent for a ring $R$:
\begin{enumerate}
       \item Every coneat left ideal of $R$ is finitely generated.
       \item Finitely $R$-injective $R$-modules are $R$-coneat injective.
       \item Finitely $R$-injective $R$-modules are self coneat injective.
       \item Absolutely self pure $R$-modules are self coneat injective.
\end{enumerate}
\end{proposition}

Now, recall that a module $M$ is called \emph{absolutely self neat} \cite{absn} if the words `finitely generated' are replaced by `maximal' in the above definition of absolutely self pure modules.
\begin{proposition} \label{maxconeat}
The following statements are equivalent for a ring $R$:
\begin{enumerate}
       \item Every maximal left ideal of $R$ is coneat in $R$.
       \item All $R$-coneat injective modules are absolutely neat.
       \item All $R$-coneat injective modules are absolutely self neat.
\end{enumerate}
\end{proposition}

A module $M$ is called \emph{principally}  (\emph{self}) injecitve  \cite{nicholson} (\cite{psi}) if every map $I \rightarrow M$ from a principal left ideal of $R$ (whose kernel is in $\Omega (M)$) can be extended to $R \rightarrow M$.
\begin{proposition} \label{coneatfg}
The following statements are equivalent for a ring $R$:
\begin{enumerate}
       \item Every coneat left ideal of $R$ is principal.
       \item Principally injective $R$-modules are $R$-coneat injective.
       \item Principally injective $R$-modules are self coneat injective.
       \item Principally self injective $R$-modules are self coneat injective.
\end{enumerate}
\end{proposition}

\begin{proposition} \label{sf}
The following are equivalent for a commutative ring $R$. 
\begin{enumerate}
       \item The ring $R$ is a von Neuman regular.
       \item Every self pure injective $R$-module is flat.
       \item Every absolutely self pure $R$-module is flat.
       \item Every absolutely self neat $R$-module is flat.
\end{enumerate}
\begin{proof} (1) $\Rightarrow$ (2) - (4) are trivial since over von Neuman regular rings, all modules are flat. (2) - (4) $\Rightarrow$ (1): Each of the conditions in (2) - (4) implies that quasi injective $R$-modules are flat. But then simple $R$-modules are flat, hence, $R$ is an SF-ring. Since $R$ is commutative, it must be von Neuman regular.
\end{proof}
\end{proposition}

\bibliographystyle{amsplain}

\end{document}